# GENERALIZATION BOUNDS FOR AVERAGED CLASSIFIERS

By Yoav Freund, Yishay Mansour[1] and Robert E. Schapire

*Columbia University, Tel-Aviv University and Princeton University*

We study a simple learning algorithm for binary classification. Instead of predicting with the best hypothesis in the hypothesis class, that is, the hypothesis that minimizes the training error, our algorithm predicts with a weighted average of all hypotheses, weighted exponentially with respect to their training error. We show that the prediction of this algorithm is much more stable than the prediction of an algorithm that predicts with the best hypothesis. By allowing the algorithm to abstain from predicting on some examples, we show that the predictions it makes when it does not abstain are very reliable. Finally, we show that the probability that the algorithm abstains is comparable to the generalization error of the best hypothesis in the class.

**1. Introduction.** Consider a binary classification learning problem. Suppose we use a hypothesis class $\mathcal{H}$ and are presented with a training set $(x_1, y_1), \ldots, (x_m, y_m)$ drawn independently from a distribution $D$ over the example domain $X \times \{-1, +1\}$. Most learning algorithms for this problem that have been studied in computational learning theory are based on identifying the hypothesis $h \in \mathcal{H}$ that minimizes the training error. One of the main problems with this approach is the phenomenon called *overfitting*. Overfitting is encountered when the hypothesis class $\mathcal{H}$ is too "large," "complex" or "flexible" relative to the size of the training set. In this case it is likely that the algorithm will find a hypothesis whose training error is very small but whose generalization error, or test error, is large. To overcome this problem, one usually uses either model-selection or regularization terms. Model selection methods try to identify the "right" complexity for $\mathcal{H}$. A regularization term is a measure of the complexity of the hypothesis $h$ that is added to the training error to define a *cost* for each hypothesis. By minimizing this cost,

Received September 2001; revised July 2003.

[1]Supported in part by a grant from the Israel Academy of Science.

the learning algorithm attempts to minimize both the training error and the amount of overfitting.

However, it is not clear that predicting with the hypothesis that minimizes the training error is indeed the only or the best prediction. One popular alternative to predicting using the single best hypothesis is to *average* the prediction of those hypotheses whose performance on the training set is close to optimal. Two popular methods of this type are Bayesian averaging [15] and bagging [4, 5]. There is considerable experimental evidence that such averaging can significantly reduce the amount of overfitting suffered by the learning algorithm. However, there is, we believe, a lack of theory for explaining this reduction.

In the context of bagging, the common explanation is based on the argument that averaging reduces the variance of the classification rule. However, as argued elsewhere [11, 18], there is currently no adequate definition of variance for classification problems. In addition, this explanation fails to take into account the effect that the complexity of the model class has on overfitting.

In the Bayesian approach the problem of overfitting is generally ignored. Instead the basic argument is that the Bayesian method is always the best method, and therefore, the only important issues are how to choose a good prior distribution and how to efficiently calculate the posterior average. However, the optimality of the Bayesian method is based on the assumption that the data we observe are *generated* according to one of the distribution models *in the chosen class of models*. While this assumption is attractive for theory, it almost never holds in practice. In practice, one usually uses relatively simple models, either because there is not enough data to estimate the "true" model, because the computational complexity is prohibitive, or because our prior knowledge of the system is only partial. Even when very complex models are used, it is rarely the case that one can assume that the data are *generated* by a model in the class. As a result, Bayesian theory is inadequate for explaining why Bayesian prediction methods are better than predicting with the best model in the class.

In this paper we propose a prediction method that is based on averaging among the empirically best classification rules. This method is similar to, but different from, the Bayesian method. The advantage of this method is that we can theoretically justify its usage without making the aforementioned Bayesian assumption that the data is generated by a distribution from a given class of distributions. Instead we make the following weaker assumptions which are common in the context of empirical error minimization methods. First, we assume that the data is generated i.i.d. according to the distribution $D$ defined above but make absolutely no assumption about $D$ other than that it is a fixed distribution. Second, we choose a class of prediction rules (mappings from the input to the binary output) and assume that there are prediction rules in that class whose probability of error (with respect to the distribution $D$) is small, *but not necessarily equal to zero.*

We deviate from the analysis used for empirical error minimization methods in our definition of a classification *rule*. In the context of a binary prediction problem, we allow the classifier *three* possible outputs. Two of them, $-1$ and $+1$, are interpreted, as before, as predictions of the label. The third, denoted by 0, should be interpreted as "no prediction" or "insufficient data."

What is the benefit of allowing the predictor this new output? The advantage is that it allows the user of the classifier to identify those examples on which overfitting might occur. For example, suppose that the best hypothesis $h^*$ in our hypothesis class $\mathcal{H}$ has an expected error of 1%. Suppose further that the size of the training set and the complexity of $\mathcal{H}$ are such that the hypothesis that minimizes the empirical error $h^*$ is likely to have a generalization error of 5%. If we use $h^*$ to make our predictions, then the most we can hope to get from a uniform-convergence type analysis is an upper bound on the generalization error that is close to 5%; we have no way of identifying *where* these errors might occur. On the other hand, if we allow the algorithm to output a zero, we can hope that the algorithm will output zero on about 4% of the input, and will be incorrect on about 1% of the data. In such a case, we say that the classifier *identifies* the locations of potential overfitting and allows the user to choose a special course of action for this case (such as referring the example back to a human to make the classification). In this case we can justifiably say that the algorithm managed to avoid overfitting. It is not misleading us into thinking that we have a classifier that is very accurate just because its error on the training set is small.

As a toy example, Figure 1 shows a tiny learning problem in which positive and negative training examples are indicated by pluses and minuses. In this example hypotheses are represented by rectangles, and we suppose that there is a large space of rectangular hypotheses, the best three of which are shown in Figure 1. Each of these makes two mistakes on this data set. However, if we take an average of hypotheses, one can imagine that it would be possible to obtain a combined classifier that abstains on all points in the shaded region where there is likely to be disagreement among the hypotheses, and predicts according to the weighted majority elsewhere. Such a combined classifier, when it does not abstain, would give nearly perfect predictions, having successfully identified the regions where errors are most likely to occur.

Of course, if the generated classifier outputs zero most of the time, then there is no benefit from having it. We need to show two things to be convinced that the addition of the new output is useful. First, we need to show that the probability of outputting a zero is of the same order as the bounds on overfitting that we would get from an analysis based on uniform convergence. Second, we need to show that when the output is $+1$ or $-1$, the probability of making a mistake is similar to the generalization error of the best hypothesis in the class. In this paper we prove that our algorithm has both these properties in the case that $\mathcal{H}$ is a finite class of models. In future work we hope to show how this work can be extended to infinite model classes.

If $\mathcal{H}$ is finite, the uniform convergence bound is the well-known Occam's razor bound [2]. If $\mathcal{H}$ is infinite, we have to resort to bounds based on VC-dimension [21]. Unfortunately, these bounds are usually very loose and provide very poor estimates for the generalization error of learning algorithms in real-world applications.

In recent years, researchers in computational learning theory have started to consider algorithms that search for a good classification rule by optimizing quantities other than the training error. Algorithms of this type include support-vector machines [21] and boosting [18] which maximize the "margin" of a linear classifier. Other work by Shawe-Taylor and Williamson [20] and McAllester [16] provide PAC-style analysis of Bayesian algorithms. Bayesian algorithms compute the posterior distribution over the space of hypotheses and predict by averaging the predictions of all hypotheses whose training error is close to the minimum. Another work that is relevant here is the work by Bousquet and Elisseeff [3] on the relationship between stability and generalization in learning classification rules.

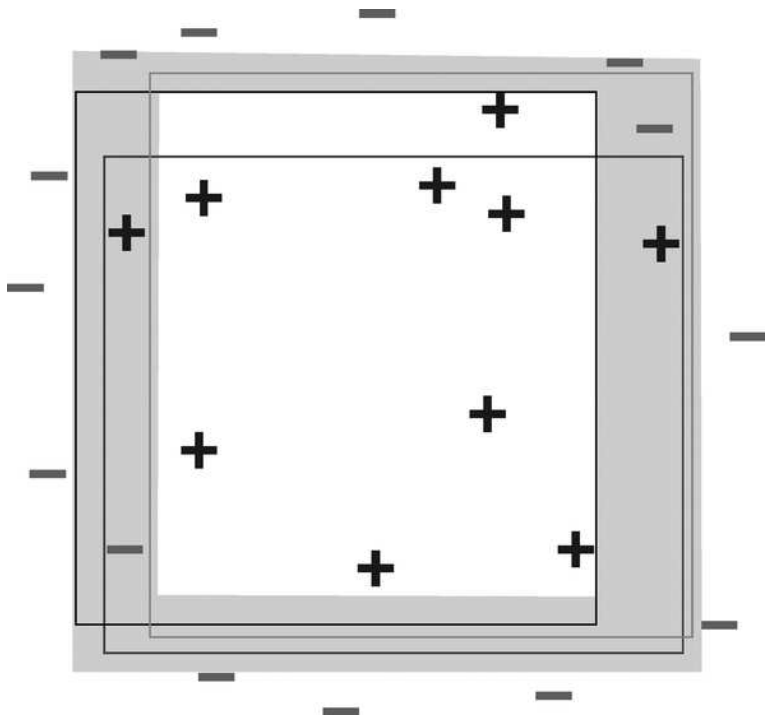

Fig. 1.  *A toy example.*

In this paper we study a learning algorithm that is very similar to the algorithm that would be suggested by Bayesian analysis but uses a slightly different formula for computing the posterior distribution. This formula is the "exponential weights" formula introduced by Littlestone and Warmuth in the context of the weighted-majority algorithm [14] and further analyzed by Cesa-Bianchi, Freund, Haussler, Helmbold, Schapire and Warmuth [6]. Note, however, that we are generating a fixed classification rule and are therefore working in the standard batch learning model and not in the online learning model.

The analysis of the algorithm consists of two parts. First, we consider, for each instance $x$, the log of the ratio of the total weight between those hypotheses that predict $+1$ on $x$ and those hypotheses that predict $-1$, where the weights depend on a parameter $\eta$. We denote this ratio by $\hat{\ell}(x)$. We prove that $\hat{\ell}(x)$ is rather insensitive to the random choice of the training set. In particular, we prove that the variation in $\hat{\ell}(x)$ is *independent of the concept class* $\mathcal{H}$! This proof is interesting because it avoids using the standard "union bound;" in fact, it altogether avoids making any uniform claim on all of the hypotheses in $\mathcal{H}$.

Using this central theorem, we can show that if $\hat{\ell}(x)$ is far from zero, then predicting with $\text{sign}(\hat{\ell}(x))$ is very stable, that is, is unlikely to change from training set to training set. More precisely, we introduce a nonstochastic quantity $\ell(x)$ and show that $\hat{\ell}(x)$ is, with high probability, very close to $\ell(x)$. Our algorithm predicts with $\text{sign}(\hat{\ell}(x))$ when $\hat{\ell}(x)$ is far from zero and abstains from prediction when $\hat{\ell}(x)$ is close to zero. We prove that the probability that this algorithm makes a prediction different from $\text{sign}(\ell(x))$ when it does not abstain is very small. On the other hand, we show that if $\mathcal{H}$ is finite and there is a hypothesis $h \in \mathcal{H}$ whose error is $\epsilon$, then we can set the parameter $\eta$ such that the error of $\text{sign}(\ell(x))$ is at most about $2\epsilon$.

The relation between our algorithm and algorithms that predict with the best hypothesis on the training set has a close correspondence to the relation between Bayesian prediction algorithms and MAP (maximum a-posteriori) algorithms. However, the analysis is carried out without making a Bayesian assumption, that is, we do not assume that the training data are generated by a model in a pre-specified class chosen by a pre-specified prior distribution. The prior and posterior distributions are internal to the algorithm and are not part of the world around it.

We hope that this paper will shed some new light on the use of algorithms that average many hypotheses such as Bayesian algorithms and averaging methods such as bagging [4, 5].

The paper is organized as follows. We start in Section 2 by describing the prediction algorithm. We give the basic analysis of the algorithm in Section 3. In Section 4 we bound the performance of $\ell(x)$ in terms of the error of the best hypothesis in the class. In Section 6 we give a bound that is uniform with respect to the learning rate parameter $\eta$ which makes it possible to choose this parameter after observing the training set. Finally, in Section 7 we outline how the ideas and results in Sections 2–4 can be extended to infinite hypothesis classes.

**2. The algorithm.** Let $D$ be a fixed but unknown distribution over $(x, y)$ pairs, where $x \in X$ and $y \in \{-1, +1\}$. Let $\mathcal{H}$ be a fixed class of hypotheses, that is, mappings from $X$ to $\{-1, +1\}$. Let $S$ denote a sample of $m$ training examples, each drawn independently at random according to $D$. We denote the *true* error of a hypothesis $h$ by $\varepsilon(h) \doteq \Pr_{(x,y) \sim D}[h(x) \neq y]$ and the estimated error according to the sample $S$ by $\hat{\varepsilon}(h) \doteq \frac{1}{m} \sum_{i=1}^{m} \mathbf{1}[h(x) \neq y]$.

The prediction algorithm that we study calculates for each hypothesis $h$ a *weight* that is defined as $w(h) \doteq e^{-\eta \hat{\varepsilon}(h)}$, where $\eta > 0$ is a parameter of the algorithm. The prediction on a new instance $x$ is defined as a function of the *empirical log ratio*:

$$\hat{\ell}_\eta(x) \doteq \frac{1}{\eta} \ln\left( \frac{\sum_{h, h(x)=+1} w(h)}{\sum_{h, h(x)=-1} w(h)} \right)$$

$$= \frac{1}{\eta} \ln\left( \frac{\sum_{h, h(x)=+1} e^{-\eta \hat{\varepsilon}(h)}}{\sum_{h, h(x)=-1} e^{-\eta \hat{\varepsilon}(h)}} \right).$$

The prediction is defined to be

$$\hat{p}_{\eta,\Delta}(x) = \begin{cases} \text{sign}(\hat{\ell}(x)), & \text{if } |\hat{\ell}(x)| > \Delta, \\ 0, & \text{otherwise,} \end{cases}$$

where $\Delta \geq 0$ is a second parameter of the algorithm. Intuitively, the parameter $\Delta$ characterizes the range of values of $\hat{\ell}_\eta(x)$ in which the training data is insufficient to make a good prediction and a better choice is to abstain. When clear from context, we generally drop the subscripts and write simply $\hat{\ell}(x)$ and $\hat{p}(x)$.

**3. Analysis of the algorithm.** For an instance $x$, we define the *true log ratio* to be

$$\ell_\eta(x) \doteq \frac{1}{\eta} \ln \frac{\sum_{h, h(x)=+1} e^{-\eta \varepsilon(h)}}{\sum_{h, h(x)=-1} e^{-\eta \varepsilon(h)}},$$

which we often write as $\ell(x)$ when $\eta$ is clear from context. The basic idea of our analysis is to show that $\hat{\ell}(x)$ must usually be close to $\ell(x)$ with high probability. In particular, we will prove the following two theorems. First, we will prove that for any fixed $x$ the difference between the empirical log ratio and the true log ratio is small:

THEOREM 1. *For any distribution $D$, any instance $x$, any $\lambda, \eta > 0$ and any $s \in \{-1, +1\}$:*

$$\Pr_{S \sim D^m}\left[s(\ell(x) - \hat{\ell}(x)) \geq 2\lambda + \frac{\eta}{8m}\right] \leq 2e^{-2\lambda^2 m}.$$

Then, in order to show that our algorithm has reasonable performance, we will transform Theorem 1 which gives a guarantee that holds with high probability for any *fixed* instance to a claim that holds with respect to a randomly chosen instance:

THEOREM 2. *For any $\delta > 0$ and $\eta > 0$, if we set*

$$\Delta = 2\sqrt{\frac{\ln(\sqrt{2}/\delta)}{m}} + \frac{\eta}{8m},$$

*then, with probability at least $1 - \delta$ over the random choice of the training set,*

$$\Pr_{(x,y) \sim D}[\hat{p}(x) \neq 0 \text{ and } \hat{p}(x) \neq \text{sign}(\ell(x))] \leq \delta.$$

This theorem guarantees that, when our algorithm predicts something different than 0 (which can be interpreted as "I do not know"), it is very likely to be making the same prediction as $\ell(x)$. Note that the statements of Theorems 1 and 2 have no dependence on the hypothesis class $\mathcal{H}$. In fact, the theorems and their proofs can be extended to infinite hypothesis classes, as discussed in Section 7.

We define some notation that will be used in the proofs. For $\mathcal{K} \subseteq \mathcal{H}$, let

$$R_\eta(\mathcal{K}) = \frac{1}{\eta} \ln\left(\sum_{h \in \mathcal{K}} e^{-\eta \varepsilon(h)}\right)$$

and let $\hat{R}_\eta(\mathcal{K})$ be the random variable

$$\hat{R}_\eta(\mathcal{K}) = \frac{1}{\eta} \ln\left(\sum_{h \in \mathcal{K}} e^{-\eta \hat{\varepsilon}(h)}\right).$$

We show that $\hat{R}_\eta(\mathcal{K})$ is close to $R_\eta(\mathcal{K})$ (with high probability) in two steps: First, we show that $\hat{R}_\eta(\mathcal{K})$ is close to its expectation $\mathrm{E}[\hat{R}_\eta(\mathcal{K})]$ with high probability. Then we show that $\mathrm{E}[\hat{R}_\eta(\mathcal{K})]$ is close to $R_\eta(\mathcal{K})$.

To prove the first result, we apply McDiarmid's theorem [17]:

THEOREM 3 (McDiarmid).   *Let $X_1, \ldots, X_m$ be independent random variables taking values in a set $V$. Let $f \colon V^m \to \mathbb{R}$ be such that, for $i = 1, \ldots, m$:*

$$|f(x_1, \ldots, x_m) - f(x_1, \ldots, x_{i-1}, x_i', x_{i+1}, \ldots, x_m)| \leq c_i$$

*for all $x_1, \ldots, x_m; x_i' \in V$. Then for $\epsilon > 0$, $s \in \{-1, +1\}$*

$$\Pr[s(f(X_1, \ldots, X_m) - \mathrm{E}[f(X_1, \ldots, X_m)]) \geq \epsilon] \leq \exp\left(-\frac{2\epsilon^2}{\sum_{i=1}^m c_i^2}\right).$$

LEMMA 1.   *Let $\mathcal{K}$ and $\hat{R}_\eta(\mathcal{K})$ be as above for a sample of size $m$. For $\eta > 0$, $\lambda > 0$ and $s \in \{-1, +1\}$,*

$$\Pr[s(\hat{R}_\eta(\mathcal{K}) - \mathrm{E}[\hat{R}_\eta(\mathcal{K})]b) \geq \lambda] \leq e^{-2\lambda^2 m}.$$

PROOF.   We apply McDiarmid's theorem with the $X_i$'s set to the labeled examples of $S$, and the function $f$ set equal to the random variable $\hat{R}_\eta(\mathcal{K})$. Let $S'$ be the sample $S$ in which one example $(x_i, y_i)$ is replaced by $(x_i', y_i')$. Let $\hat{\varepsilon}'(h)$ be the empirical error of $h$ on $S'$, and let

$$\hat{R}_\eta'(\mathcal{K}) = \frac{1}{\eta} \ln\left(\sum_{h \in \mathcal{K}} e^{-\eta \hat{\varepsilon}'(h)}\right).$$

Then

$$\begin{aligned}
\hat{R}_\eta'(\mathcal{K}) - \hat{R}_\eta(\mathcal{K}) &= \frac{1}{\eta} \ln\left(\frac{\sum_{h \in \mathcal{K}} e^{-\eta \hat{\varepsilon}'(h)}}{\sum_{h \in \mathcal{K}} e^{-\eta \hat{\varepsilon}(h)}}\right) \\
&\leq \frac{1}{\eta} \ln\left(\max_{h \in \mathcal{K}} e^{-\eta(\hat{\varepsilon}'(h) - \hat{\varepsilon}(h))}\right) \\
&= \max_{h \in \mathcal{K}}(\hat{\varepsilon}'(h) - \hat{\varepsilon}(h)) \leq \frac{1}{m}.
\end{aligned}$$

The first inequality uses the fact that $(\sum_i a_i)/(\sum_i b_i) \leq \max_i a_i / b_i$ for positive $a_i$'s and $b_i$'s. The second inequality uses the fact that changing one example can change the empirical error by at most $1/m$.

By the symmetry of this argument, $|\hat{R}_\eta(\mathcal{K}) - \hat{R}_\eta'(\mathcal{K})| \leq 1/m$. Plugging in $c_i = 1/m$ in McDiarmid's theorem gives the result.   □

LEMMA 2.   *Let $\mathcal{K}$, $R_\eta(\mathcal{K})$ and $\hat{R}_\eta(\mathcal{K})$ be as above for a sample of size $m$. Then for $\eta > 0$,*

$$R_\eta(\mathcal{K}) \leq \mathrm{E}[\hat{R}_\eta(\mathcal{K})] \leq R_\eta(\mathcal{K}) + \frac{\eta}{8m}.$$

PROOF. For the lower bound on $\mathrm{E}[\hat{R}_\eta(\mathcal{K})]$, let $\mathcal{K} = \{h_1, \ldots, h_N\}$. For $\mathbf{x} \in \mathbb{R}^N$, let

$$g(\mathbf{x}) = \ln\left(\sum_{i=1}^N e^{x_i}\right).$$

Then $g$ is convex: Given $\alpha \in (0,1)$ and $\mathbf{x}, \mathbf{y} \in \mathbb{R}^N$, let $p = 1/\alpha$, $q = 1/(1-\alpha)$, and define $r_i = e^{\alpha x_i}$ and $s_i = e^{(1-\alpha)y_i}$. Since $1/p + 1/q = 1$, by Hölder's inequality,

$$\sum_i r_i s_i \le \left(\sum_i r_i^p\right)^{1/p} \left(\sum_i s_i^q\right)^{1/q}.$$

Plugging in definitions and taking logarithms, this is equivalent to

$$g(\alpha \mathbf{x} + (1-\alpha)\mathbf{y}) \le \alpha g(\mathbf{x}) + (1-\alpha)g(\mathbf{y}),$$

so $g$ is convex as claimed.

Therefore, by Jensen's inequality,

$$\begin{aligned}
\eta \mathrm{E}[\hat{R}_\eta(\mathcal{K})]\mathcal{K} &= \mathrm{E}[g(\langle -\eta \hat{\varepsilon}(h_1), \ldots, -\eta \hat{\varepsilon}(h_N)\rangle)] \\
&\ge g(\langle -\eta \mathrm{E}[\hat{\varepsilon}(h_1)], \ldots, -\eta \mathrm{E}[\hat{\varepsilon}(h_N)]\rangle) \\
&= g(\langle -\eta \varepsilon(h_1), \ldots, -\eta \varepsilon(h_N)\rangle) = \eta R_\eta(\mathcal{K}).
\end{aligned}$$

To prove the upper bound on $\mathrm{E}[\hat{R}_\eta(\mathcal{K})]$, we have by Jensen's inequality (applied to the concave log function),

$$\begin{aligned}
(1) \qquad \mathrm{E}[\hat{R}_\eta(\mathcal{K})] &= \frac{1}{\eta} \mathrm{E}\left[\ln\left(\sum_{h \in \mathcal{K}} e^{-\eta \hat{\varepsilon}(h)}\right)\right] \\
&\le \frac{1}{\eta} \ln\left(\sum_{h \in \mathcal{K}} \mathrm{E}[e^{-\eta \hat{\varepsilon}(h)}]\right).
\end{aligned}$$

Fix $h$ and let $\varepsilon = \varepsilon(h)$ and $\hat{\varepsilon} = \hat{\varepsilon}(h)$. Let $Z_i$ be a Bernoulli random variable that is 1 if $h(x_i) \ne y_i$ and 0 otherwise. Then we can write

$$\begin{aligned}
\mathrm{E}[e^{\eta(\varepsilon - \hat{\varepsilon})}] &= \mathrm{E}\left[\exp\left(\frac{\eta}{m}\sum_{i=1}^m (\varepsilon - Z_i)\right)\right] \\
&= \prod_{i=1}^m \mathrm{E}\left[\exp\left(\frac{\eta}{m}(\varepsilon - Z_i)\right)\right] \\
&\le (e^{\eta^2/8m^2})^m = e^{\eta^2/8m}.
\end{aligned}$$

The second equality uses independence of the $Z_i$'s. The last step uses the fact, proved by Hoeffding [13], that for any random variable $X$ with $\mathrm{E}[X] = 0$ and $a \le X \le b$,

$$\mathrm{E}[e^X] \le e^{(b-a)^2/8}.$$

Here we let $X = (\eta/m)(\varepsilon - Z_i)$.

Thus, $\mathrm{E}[e^{-\eta \hat{\varepsilon}(h)}] \le e^{\eta^2/8m} e^{-\eta \varepsilon(h)}$. Combined with (1), this gives that

$$\mathrm{E}[\hat{R}_\eta(\mathcal{K})] \mathcal{K} \le \frac{1}{\eta} \ln \left( e^{\eta^2/8m} \sum_{h \in \mathcal{K}} e^{-\eta \varepsilon(h)} \right) = R_\eta(\mathcal{K}) + \frac{\eta}{8m}$$

as claimed.  □

PROOF OF THEOREM 1.   Given $x$, we partition the hypothesis set $\mathcal{H}$ into two. The subset $\mathcal{K}$ includes the hypotheses $h$ such that $h(x) = +1$ and its complement $\mathcal{K}^c$ includes all $h$ for which $h(x) = -1$. We can now write

$$
\begin{aligned}
(2) \quad \ell(x) - \hat{\ell}(x) &= \frac{1}{\eta} \ln \left( \frac{\sum_{h \in \mathcal{K}} e^{-\eta \varepsilon(h)}}{\sum_{h \in \mathcal{K}} e^{-\eta \hat{\varepsilon}(h)}} \right) + \frac{1}{\eta} \ln \left( \frac{\sum_{h \in \mathcal{K}^c} e^{-\eta \hat{\varepsilon}(h)}}{\sum_{h \in \mathcal{K}^c} e^{-\eta \varepsilon(h)}} \right) \\
&= R_\eta(\mathcal{K}) - R_\eta(\mathcal{K}^c) - \hat{R}_\eta(\mathcal{K}) + \hat{R}_\eta(\mathcal{K}^c).
\end{aligned}
$$

Combining Lemmas 1 and 2, we find that

$$(3) \qquad \Pr[R_\eta(\mathcal{K}) - \hat{R}_\eta(\mathcal{K}) > \lambda] \le e^{-2\lambda^2 m}$$

and

$$(4) \qquad \Pr\left[ \hat{R}_\eta(\mathcal{K}^c) - R_\eta(\mathcal{K}^c) > \lambda + \frac{\eta}{8m} \right] \le e^{-2\lambda^2 m}.$$

Combining (2)–(4), we prove the claim for $s = +1$. The proof for $s = -1$ is almost identical.  □

LEMMA 3.   *For any distribution $D$, any $\lambda, \eta > 0$ and any $s \in \{-1, +1\}$, the probability over samples $S \sim D^m$ that*

$$\Pr_{(x,y) \sim D} \left[ s(\ell(x) - \hat{\ell}(x)) \ge 2\lambda + \frac{\eta}{8m} \right] \ge \sqrt{2} e^{-\lambda^2 m}$$

*is at most $\sqrt{2} e^{-\lambda^2 m}$.*

PROOF.   Since Theorem 1 holds for all $x$, it also holds for a random $x$. Thus,

$$
\begin{aligned}
\mathop{\mathrm{E}}_{S \sim D^m} &\left[ \Pr_{(x,y) \sim D} \left[ s(\ell(x) - \hat{\ell}(x)) \ge 2\lambda + \frac{\eta}{8m} \right] \right] \\
&= \mathop{\mathrm{E}}_{(x,y) \sim D} \left[ \Pr_{S \sim D^m} \left[ s(\ell(x) - \hat{\ell}(x)) \ge 2\lambda + \frac{\eta}{8m} \right] \right] \\
&\le 2 e^{-2\lambda^2 m}.
\end{aligned}
$$

The lemma now follows using Markov's inequality. □

Theorem 2 follows immediately from Lemma 3.

**4. Performance relative to the best hypothesis.** We now show that there exists a setting of $\eta$ and $\Delta$ that yields performance guarantees relative to the best hypothesis in the class. We compare these guarantees to those given by the Occam argument [2] for the algorithm that uses a hypothesis that minimizes the empirical error rate.

In Lemma 3 we showed that the value of $\hat{\ell}(x)$ is, with high probability, close to $\ell(x)$. We now show that, with respect to the *actual* distribution $D$, the sign of $\ell(x)$ is closely related to that of the best hypothesis in $\mathcal{H}$. By combining these theorems, we show that the generalization error of our algorithm is close to that of the best hypothesis in $\mathcal{H}$.

Note that the following theorem does not involve the training set in any way; it is a claim about $y\ell(x)$ which is a deterministic function of $(x, y)$. Intuitively, for large enough values of $\eta$, the function $\ell(x)$ essentially averages the best hypotheses from $\mathcal{H}$. In the worse case, as we show in Section 5, this can at most double the error. The following theorem gives a detailed tradeoff between all the parameters.

THEOREM 4. *Let $\mathcal{H}$ be a finite hypothesis class and let $\epsilon$ be the error of the best hypothesis in $\mathcal{H}$ with respect to the distribution $D$ over the examples, that is, $\epsilon = \min\{\varepsilon(h) : h \in \mathcal{H}\}$. Let $\eta > 0$ and $\Delta \geq 0$ be such that $\Delta\eta \leq 1/2$. Then for any $\gamma \geq \ln(8|\mathcal{H}|)/\eta$,*

$$\Pr_{(x,y)\sim D}[y\ell(x) \leq 0] \leq 2(1 + 2|\mathcal{H}|e^{-\eta\gamma})(\epsilon + \gamma),$$

*and*

$$\Pr_{(x,y)\sim D}[y\ell(x) \leq 2\Delta] \leq (1 + e^{2\Delta\eta})(1 + 2|\mathcal{H}|e^{\eta(2\Delta-\gamma)})(\epsilon + \gamma)$$

$$\leq 4(1 + 2|\mathcal{H}|e^{\eta(2\Delta-\gamma)})(\epsilon + \gamma).$$

PROOF. We partition the hypotheses in $\mathcal{H}$ into two sets according to their true error. We call those hypotheses whose error is smaller than $\epsilon + \gamma$ *strong* and the other hypotheses *weak*.

We denote by $W_w$ the total weight of the weak hypotheses:

$$W_w = \frac{1}{Z} \sum_{h \in \mathcal{H} : \varepsilon(h) \geq \epsilon+\gamma} e^{-\eta\varepsilon(h)},$$

where

$$Z = \sum_{h \in \mathcal{H}} e^{-\eta\varepsilon(h)}.$$

To upper bound $W_w$, note that we always have at least one strong hypothesis, namely, the one that achieves $\varepsilon(h) = \epsilon$. Thus,

$$(5) \qquad W_w \leq \frac{|\mathcal{H}|e^{-\eta(\epsilon+\gamma)}}{e^{-\eta\epsilon}} = |\mathcal{H}|e^{-\eta\gamma}.$$

From the assumption that $\gamma \geq \ln(8|\mathcal{H}|)/\eta$, we get that $W_w \leq 1/8$.

For a given example $(x, y)$, we partition the strong hypotheses into two subsets according to whether or not the hypothesis gives the correct prediction on $(x, y)$. We denote the total weight of these subsets by

$$W_s^+(x, y) = \frac{1}{Z} \sum_{h \in \mathcal{H} \,:\, \varepsilon(h) < \epsilon + \gamma, h(x) = y} e^{-\eta\varepsilon(h)},$$

$$W_s^-(x, y) = \frac{1}{Z} \sum_{h \in \mathcal{H} \,:\, \varepsilon(h) < \epsilon + \gamma, h(x) \neq y} e^{-\eta\varepsilon(h)}.$$

By the definition of $Z$, for any $(x, y)$,

$$W_s^+(x, y) + W_s^-(x, y) + W_w = 1.$$

We now prove the second part of the theorem; the first part follows from the second part by setting $\Delta = 0$. We first bound $y\ell(x)$ using $W_w$, $W_s^+(x, y)$ and $W_s^-(x, y)$:

$$y\ell(x) \geq \frac{1}{\eta} \ln\left(\frac{W_s^+(x, y)}{W_s^-(x, y) + W_w}\right).$$

Thus, $y\ell(x) \leq 2\Delta$ implies

$$\frac{W_s^-(x, y) + W_w}{1 - (W_s^-(x, y) + W_w)} \geq e^{-2\Delta\eta}$$

or, equivalently,

$$W_s^-(x, y) + W_w \geq \frac{1}{1 + e^{2\Delta\eta}} \doteq c.$$

We denote by $h \sim \mathcal{S}$ the random choice of a hypothesis from the strong set with probability $e^{-\eta\varepsilon(h)}/Z_s$, where $Z_s$ normalizes the weights *within* the strong set to sum to 1. We find that

$$\Pr_{(x,y) \sim D}[y\ell(x) \leq 2\Delta] \leq \Pr_{(x,y) \sim D}\left[\frac{W_s^-(x, y)}{W_s^-(x, y) + W_s^+(x, y)} \geq \frac{c - W_w}{1 - W_w}\right]$$

$$= \Pr_{(x,y) \sim D}\left[\Pr_{h \sim S}[h(x) \neq y] \geq \frac{c - W_w}{1 - W_w}\right]$$

$$(6) \qquad \leq \mathop{\mathrm{E}}_{(x,y) \sim D}\left[\Pr_{h \sim S}[h(x) \neq y]\right]\frac{1 - W_w}{c - W_w}$$

$$(7) \qquad = \underset{h \sim S}{\mathrm{E}} \left[ \underset{(x,y) \sim D}{\mathrm{Pr}} [h(x) \neq y] \right] \frac{1 - W_w}{c - W_w}$$

$$(8) \qquad \leq (\epsilon + \gamma) \frac{1 - W_w}{c - W_w}$$

$$(9) \qquad \leq (\epsilon + \gamma)(1 + e^{2\Delta\eta})(1 + 2W_w e^{2\Delta\eta}).$$

Equations (6) and (7) use Markov's inequality and Fubini's theorem. Equation (8) follows from the fact that $\varepsilon(h) < \epsilon + \gamma$ for every strong hypothesis. Equation (9) uses our assumptions that $\Delta\eta \leq 1/2$ and $W_w \leq 1/8$ together with the inequality $(1 - x)/(1 - x(1 + r)) \leq 1 + 2xr$ for $x > 0$, $r > 0$ and $x(1 + r) \leq 1/2$ (with $x = W_w$ and $r = e^{2\Delta\eta}$).

Combining this bound with (5) proves the second statement of the theorem. □

**5. Discussion.** We now discuss the implications of Theorems 2 and 4. We start with a corollary of Theorem 4 for a specific setting of the parameters $\eta$ and $\Delta$ as a function of the sample size $m$, the size of the hypothesis class $\mathcal{H}$ and the reliability parameter $\delta$.

COROLLARY 1. *Let $1/2 > \theta > 0$, $\delta > 0$ and*

$$\eta = \ln(8|\mathcal{H}|)m^{1/2-\theta}; \qquad \Delta = 2\sqrt{\frac{\ln(\sqrt{2}/\delta)}{m}} + \frac{\ln(8|\mathcal{H}|)}{8m^{1/2+\theta}}.$$

*For $m \geq 8$,*

$$\underset{(x,y) \sim D}{\mathrm{Pr}}[y\ell(x) \leq 0] \leq \left(2 + \frac{1}{4m}\right) \left(\epsilon + \frac{1}{m^{1/2-\theta}} + \frac{\ln m}{m^{1/2-\theta} \ln 8|\mathcal{H}|}\right),$$

*and for*

$$m \geq \left[8\sqrt{\ln\left(\frac{\sqrt{2}}{\delta}\right)} \ln(8|\mathcal{H}|)\right]^{1/\theta},$$

*we have*

$$\underset{(x,y) \sim D}{\mathrm{Pr}}[y\ell(x) \leq 2\Delta] \leq 5\left(\epsilon + 2\Delta + \frac{1}{m^{1/2-\theta}}\right).$$

PROOF. To prove the corollary, we use Theorem 4 with two different settings of $\gamma$. The first bound is a result of choosing $\gamma = 1/m^{1/2-\theta} + \ln m/(m^{1/2-\theta} \ln 8|\mathcal{H}|)$, and the second is a result of choosing $\gamma = 2\Delta + m^{\theta-1/2}$. □

We now discuss the significance of each statement in the corollary. Let us fix the reliability parameter $\delta$.

The first statement of Corollary 1 shows that the sign of the true log ratio is a reasonably good proxy for the best hypothesis in the class, denoted $h^*$. Specifically, the error of $\text{sign}(\ell(x))$ is

$$2\varepsilon(h^*) + O\left(\frac{\ln(m)}{m^{1/2-\theta}}\right).$$

Let us separate between abstaining and making a mistake. If the algorithm outputs 0 we say that it "abstained," while if it outputs $-1$ or $+1$ and this label does not agree with the actual label of the example, then we say that it "made a mistake." Combining this with the statement of Theorem 2, we find that the probability that our algorithm makes a mistake on a test example is bounded by

$$(10) \qquad 2\varepsilon(h^*) + O\left(\frac{\ln(m)}{m^{1/2-\theta}}\right) + \delta.$$

Note that this bound is *independent* of $|\mathcal{H}|$.

In comparison, the upper bound on the hypothesis that minimizes the empirical risk is

$$(11) \qquad \varepsilon(h^*) + O\left(\sqrt{\frac{\ln(|\mathcal{H}|/\delta)}{m}}\right).$$

We see that the dependence on $m$ here is slightly better, but the bound depends on the hypothesis class, which is what we expect from an algorithm that cannot abstain.

For our algorithm, the dependence on $|\mathcal{H}|$ instead appears in the bound on the probability of abstaining on a test example; this is given in the second statement of the corollary. Combining that statement with Lemma 3, we find that for

$$m = \Omega((\sqrt{\ln(1/\delta)}\ln(|\mathcal{H}|))^{1/\theta})$$

our algorithm will predict zero with probability at most

$$5\varepsilon(h^*) + O\left(\frac{\sqrt{\ln(1/\delta)} + \ln(|\mathcal{H}|)}{m^{1/2-\theta}}\right).$$

This bound is similar to the Occam bound (11), but the choice of $\theta$ makes an important difference in the dependence on $m$.

In effect, we are replacing one type of guarantee with a different one. In the traditional analysis that is based on uniform convergence theory, the guarantee is of the form "the error of the classification rule is at most $\epsilon + O(\ln(|H|)/\sqrt{m})$." Our algorithm is one for which there are two guarantees. First, we can say that "the error of the classification rule, when this rule makes a nonzero prediction, is at most $2\epsilon + \tilde{O}(1/\sqrt{m})$ (no dependence on the size of $H$ here). Second, we can show that the probability that the classification rule will generate a 0 ("I do not know" prediction) is upper bounded by $5\epsilon + \tilde{O}(\ln(|H|)/\sqrt{m})$. This second bound does depend on the size of $H$. Note that this quantity (the probability of predicting 0) can be estimated from an *unlabeled* set of instances. Unlike the event of a classification mistake, which depends both on the predicted label and the actual label, the event of predicting 0 does not depend on the actual label. In practice, unlabeled data is usually much more plentiful than labeled data. Therefore, in practice, we can estimate the probability of abstaining directly and do not need to use a priori bounds.

We now argue that the factor of 2 in front of the error of the best hypothesis in the class which appears in the first part of the corollary is necessary. Suppose that the input domain $X$ is partitioned into two parts $A_1$ and $A_2$, such that $D(A_1) = 1 - 2\epsilon$ and $D(A_2) = 2\epsilon$. Suppose that all the hypotheses in $\mathcal{H}$ predict correctly on instances in $A_1$. For each $x \in A_2$, the prediction of each hypothesis is chosen independently at random to be correct with probability $1/2 - \eta\Delta$ and incorrect with probability $1/2 + \eta\Delta$. (Suppose further the number of elements in $A_2$ and the number of hypotheses are sufficiently large so that on most of the points in $A_2$ the actual fraction of correct predictions is sufficiently close to $1/2 - \eta\Delta$.) In this case each of the hypotheses in $\mathcal{H}$ has error close to $2\epsilon(1/2 + \eta\Delta) \approx \epsilon(1 + O(m^{-\theta}))$. This also implies that all of the hypotheses have approximately the same weight.

Consider now the value of $\ell(x)$ for $x \in A_2$. As the weights of all of the hypotheses are similar, we get that

$$\forall\, x \in A_2, \qquad y\ell(x) \approx \frac{1}{\eta}\ln\!\left(\frac{1/2 - \eta\Delta}{1/2 + \eta\Delta}\right) \approx -4\Delta.$$

As $\hat{\ell}$ is likely to be very close to $\ell$, we conclude that for $x \in A_2$ our algorithm will usually make a nonzero prediction that is incorrect. In other words, our algorithm will have a prediction error of about $2\epsilon$ while each of the hypotheses has error of about $\epsilon$.

It may seem impossible that the bound in (10) is independent of the number of hypotheses. First, one should recall that a similar phenomenon exists in the large margin analysis for hyperplanes, where the generalization error depends only on the margin and not on the dimension of the class. One should not interpret our result as suggesting that overfitting can never happen, regardless of the complexity of the hypothesis space. In truth, if the hypothesis space is too complex, the algorithm will simply abstain more often. For example, suppose that the hypothesis space consists of *all* binary functions on a finite domain. For any set of training examples, there is a function that has zero training error (assuming no example appears twice with different labels). However, we expect any algorithm to be unable to predict the label of a new test example. Indeed, in this case, our algorithm will abstain on all unseen examples [since $\hat{\ell}(x)$ is exactly zero outside the training set].

Using the size of the hypothesis class as the measure of its complexity is clearly a very rough upper bound. For example, consider the case in which a large fraction of the hypotheses in $H$ are all equal, or almost equal, to a single function $h^*$. It is not hard to see that in this case our prediction algorithm, as stated, will have a strong bias towards predicting like $h^*$. This bias can be removed by replacing the set of almost identical hypotheses by the single hypothesis $h^*$. Doing this also improves the guaranteed performance bounds because it reduces $|H|$. A systematic way for removing this type of bias is to replace $H$ with an $\epsilon$-*net* that covers it. In other words, find a set of functions $H_\epsilon$ such that for any $h \in H$ there exists $f \in H_\epsilon$ such that $\Pr_{(x,y) \sim D}[h(x) \neq f(x)] \leq \epsilon$. Of course, choosing an $\epsilon$-cover requires knowledge of the marginal distribution over $x$ defined by $D$ and is a nontrivial computational problem. Potential future research regarding the use of $\epsilon$-covers in conjunction with our prediction algorithm is discussed in Section 9.

Finally, Theorem 4 shows that the error of our predictor cannot be much worse than twice the error of the best hypothesis. On the other hand, it is possible in some favorable situations for our predictor to significantly outperform the best hypothesis. For example, suppose that there is an $h^* \in \mathcal{H}$ such that $\varepsilon(h^*) = 1/8$, and that for each $h \in \mathcal{H}' = \mathcal{H} - \{h^*\}$, we have $\varepsilon(h) = 1/4$. Suppose further that for each $x$, the fraction of $h \in \mathcal{H}'$ with the right label is $3/4$. Choosing the hypothesis with lowest observed error would give, hopefully, the hypothesis $h^*$ that has an error rate of $1/8$. In our setting, for a labeled example $(x, y)$, if $h^*(x) = y$, then

$$y\ell(x) = \frac{1}{\eta} \ln \left( \frac{e^{-\eta/8} + (3/4)|\mathcal{H}'|e^{-\eta/4}}{(1/4)|\mathcal{H}'|e^{-\eta/4}} \right)$$

$$= \frac{1}{\eta} \ln \left( 3 + \frac{4e^{\eta/8}}{|\mathcal{H}'|} \right).$$

Thus, for $\eta = 1$, we have $y\ell(x) = \ln(3 + 4e^{1/8}/|\mathcal{H}'|)$. Similarly, if $h^*(x) \neq y$, we have $y\ell(x) \geq \ln(3 - 12e^{1/8}/|\mathcal{H}'|)$. Note that this implies that $p_{1,0}(x)$ correctly classifies all the examples (for $|\mathcal{H}|$ large). Theorem 1, with $\lambda$ set to a constant, then guarantees for $m = O(\lg 1/\delta)$ that $\hat{p}_{1,0}(x)$ has an error rate of at most $\delta$. The important point here is that by averaging a large number of suboptimal hypotheses we achieve a prediction accuracy that is better than that of the optimal single hypothesis $h^*$.

An interesting question was raised by one of the reviewers: why are we comparing the performance of our algorithm to that of the optimal single prediction rule when, in fact, one would expect a rule that is a combination of many prediction rules to perform much better than any single rule? Our answer is that in this work we wanted to relate our bounds to those that are proven using uniform-convergence analysis of the type advocated by Vapnik [21], and those have as their "gold standard" the performance of the optimal hypothesis. A natural direction for future research would be to compare the performance of our algorithm to that of the rule:

$$\mathrm{sign}\left(\lim_{\eta \to \infty} \ell_\eta(x)\right),$$

which is the analog of our prediction rule when the distribution is known (or equivalently, in the limit of an infinite number of training examples). However, it is not clear whether this is the correct gold standard to use.

**6. Uniform bounds.** The bound given in Lemma 1 applies to the case in which the parameter $\eta$ is fixed ahead of time so that $\hat{R}_\eta(\mathcal{K})$ converges to $\mathrm{E}[\hat{R}_\eta(\mathcal{K})]$ for only a single value of $\eta$. In the next lemma we show that on a single sample this convergence is likely to take place for *all* values of $\eta \geq 1$ simultaneously. (We can prove a similar result for $\eta > 0$ using a slightly more complicated proof. However, because $\eta$ is typically large in this paper, we omit this proof.) The proof of this is primarily taken from Allwein, Schapire and Singer [1].

LEMMA 4.    *Let* $\mathcal{K}$ *and* $\hat{R}_\eta(\mathcal{K})$ *be as above for a sample of size* $m$. *For* $\lambda > 0$,

$$\Pr[\exists\, \eta \geq 1 : |\hat{R}_\eta(\mathcal{K}) - \mathrm{E}[\hat{R}_\eta(\mathcal{K})]| \geq \lambda] \leq \frac{8 \ln |\mathcal{K}|}{\lambda} e^{-\lambda^2 m/2}.$$

The proof is given in the Appendix.

We can now state the following theorems similar to Theorems 1 and 2. These theorems show that it is possible to design an algorithm that chooses $\eta$ *after* the sample has been chosen without paying a large penalty in accuracy.

THEOREM 5.    *Let* $\mathcal{K}$ *and* $\hat{R}_\eta(\mathcal{K})$ *be as above for a sample of size* $m$. *For any distribution* $D$, *any* $\lambda > 0$ *and any* $s \in \{-1, +1\}$,

$$\Pr_{S \sim D^m}\left[\exists\, \eta \geq 1 : s(\ell_\eta(x) - \hat{\ell}_\eta(x)) \geq 2\lambda + \frac{\eta}{8m}\right] \leq \frac{8 \ln |\mathcal{K}|}{\lambda} e^{-\lambda^2 m/2}.$$

THEOREM 6.    *For any* $\delta > 0$, *if we set*

$$\Delta_\eta = 2\sqrt{\frac{2}{m}\ln\left(\frac{16m \ln |\mathcal{H}|}{\delta^2}\right)} + \frac{\eta}{m},$$

*then, with probability at least $1 - \delta$ over the random choice of the training set, for all $\eta \geq 1$*

$$\Pr_{(x,y) \sim D}[\hat{p}_{\eta, \Delta_\eta}(x) \neq 0 \text{ and } \hat{p}_{\eta, \Delta_\eta}(x) \neq \text{sign}(\ell_\eta(x))] \leq \delta.$$

## 7. Infinite hypothesis classes.

The ideas and results of Sections 2–4 can be directly extended to infinite, even uncountable, hypothesis spaces. To make this extension, we need to add as a parameter of the algorithm a finite measure $\mu$ over the hypothesis space $\mathcal{H}$. For convenience, we assume in fact that $\mu$ is a probability measure so that

$$\mu(\mathcal{H}) = \int_{\mathcal{H}} d\mu = 1.$$

Naturally, we will require certain measurability assumptions so that everything is measurable that needs to be so. For our purposes, it is sufficient that the following sets are measurable:

$$\{h \in \mathcal{H} : h(x) = +1\}, \qquad \text{for all } x \in X,$$
$$\{h \in \mathcal{H} : \varepsilon(h) < \epsilon\}, \qquad \text{for all } \epsilon \in \mathbb{R}.$$

In other words, these sets are assumed to be elements of the $\sigma$-algebra over which the measure $\mu$ is defined.

The results for finite $\mathcal{H}$ presented earlier in the paper are, of course, a special case in which $\mu$ is the uniform discrete measure $\mu(\mathcal{K}) = |\mathcal{K}|/|\mathcal{H}|$ for all $\mathcal{K} \subseteq \mathcal{H}$.

Formally, the measure $\mu$ is used much like a Bayesian prior. However, unlike a prior, we do *not* assume that there is a target hypothesis in $\mathcal{H}$ that has been chosen randomly according to $\mu$.

The algorithm in Section 2 can now be extended by simply redefining the empirical log ratio to be

$$\hat{\ell}_\eta(x) \doteq \frac{1}{\eta} \ln \left( \frac{\int_{\{h : h(x) = +1\}} w(h) \, d\mu}{\int_{\{h : h(x) = -1\}} w(h) \, d\mu} \right),$$

where as usual $w(h) \doteq e^{-\eta \hat{\varepsilon}(h)}$ and the integral is the Lebesgue integral with regard to the probability measure. The true log ratio $\ell_\eta(x)$ is redefined analogously.

To prove Theorems 1 and 2 and Lemmas 1–3 in this more general setting, we simply need to replace each sum of the form $\sum_{h \in \mathcal{K}} f(h)$ by the integral $\int_{\mathcal{K}} f(h) \, d\mu$ for measurable sets $\mathcal{K}$. [If $\mathcal{K}$ has measure zero, then $R_\eta(\mathcal{K})$ and $\hat{R}_\eta(\mathcal{K})$ are both defined to be zero.]

The only potential difficulty occurs in proving in Lemma 2 that $R_\eta(\mathcal{K}) \leq \mathrm{E}[\hat{R}_\eta(\mathcal{K})]$. When $\mathcal{K}$ is finite, we can simply apply Jensen's inequality to a function of $|\mathcal{K}|$ real variables. When $\mathcal{K}$ is infinite, however, this may be a problem since standard forms of Jensen's inequality do not apply. Nevertheless, we can effectively reduce to the finite case as follows:

Let $\delta > 0$. Let

$$\mathcal{B}_i = \{h \in \mathcal{K} : i\delta \le \varepsilon(h) < (i+1)\delta\}.$$

Since $\varepsilon(h) \in [0,1]$, $B_0, \dots, B_k$ form a partition of $\mathcal{K}$ for $k = \lfloor 1/\delta \rfloor$. For $\mu(\mathcal{B}_i) > 0$, define $\tilde{\varepsilon}_i$ to be a random variable that is the average of $\hat{\varepsilon}(h)$ over $h \in \mathcal{B}_i$, that is,

$$\tilde{\varepsilon}_i \doteq \frac{\int_{\mathcal{B}_i} \hat{\varepsilon}(h)\,d\mu}{\mu(\mathcal{B}_i)}.$$

Then

$$\mathrm{E}\left[\tilde{\varepsilon}_i\right] = \frac{\int_{\mathcal{B}_i} \varepsilon(h)\,d\mu}{\mu(\mathcal{B}_i)} \le (i+1)\delta.$$

Combined with the fact that $\varepsilon(h) \ge i\delta$, for $h \in \mathcal{B}_i$, gives

$$\begin{aligned}
\int_{\mathcal{K}} e^{-\eta\varepsilon(h)}\,d\mu &\le \sum \mu(\mathcal{B}_i) e^{-\eta i\delta} \\
&\le \sum \mu(\mathcal{B}_i) e^{-\eta(\mathrm{E}[\tilde{\varepsilon}_i] - \delta)} \\
&= e^{\eta\delta} \sum \mu(\mathcal{B}_i) e^{-\eta\,\mathrm{E}[\tilde{\varepsilon}_i]},
\end{aligned}$$

where it is understood that all sums are over $i$ for which $\mu(\mathcal{B}_i) > 0$. Thus,

$$\begin{aligned}
R_\eta(\mathcal{K}) &= \frac{1}{\eta} \ln \int_{\mathcal{K}} e^{-\eta\varepsilon(h)}\,d\mu \\
&\le \delta + \frac{1}{\eta} \ln \sum \mu(\mathcal{B}_i) e^{-\eta\,\mathrm{E}[\tilde{\varepsilon}_i]} \\
(12) \qquad &\le \delta + \frac{1}{\eta} \mathrm{E}\left[\ln \sum \mu(\mathcal{B}_i) e^{-\eta\tilde{\varepsilon}_i}\right] \\
&= \delta + \frac{1}{\eta} \mathrm{E}\left[\ln \sum \mu(\mathcal{B}_i) \exp\left(-\eta \frac{\int_{\mathcal{B}_i} \hat{\varepsilon}(h)\,d\mu}{\mu(\mathcal{B}_i)}\right)\right] \\
(13) \qquad &\le \delta + \frac{1}{\eta} \mathrm{E}\left(\ln \sum \mu(\mathcal{B}_i) \frac{\int_{\mathcal{B}_i} e^{-\eta\hat{\varepsilon}(h)\,d\mu}}{\mu(\mathcal{B}_i)}\right) \\
&= \delta + \frac{1}{\eta} \mathrm{E}\left(\ln \int_{\mathcal{K}} e^{-\eta\hat{\varepsilon}(h)\,d\mu}\right) \\
&= \delta + \mathrm{E}[\hat{R}_\eta(\mathcal{K})].
\end{aligned}$$

Equation (12) uses Jensen's inequality applied to the convex function

$$\mathbf{x} \mapsto \ln \sum_i \mu(\mathcal{B}_i) e^{x_i}.$$

(Convexity follows from a minor modification of the proof given in Lemma 2 for the function $g$.) Equation (13) applies Jensen's inequality to the convex function $e^x$. Since $\delta$ is arbitrary, the result follows.

The results in Section 4 compare performance to that of the best single hypothesis. When $\mathcal{H}$ is infinite, this comparison may be meaningless since this single hypothesis is likely to have measure zero. Moreover, the bounds in Section 4 are in terms of $|\mathcal{H}|$ which will now be infinite.

Therefore, rather than comparing to a single best hypothesis, we compare to a set of good hypotheses. In particular, for any $\epsilon > 0$, let $V_\epsilon$ be the volume of all hypotheses with error at most $\epsilon$:

$$V_\epsilon \doteq \mu(\{h : \varepsilon(h) \leq \epsilon\}).$$

Then throughout this section we need to replace $|\mathcal{K}|$ with $1/V_\epsilon$.

Specifically, the generalization of Theorem 4 becomes the following:

THEOREM 7. *Let $\mathcal{H}$ be any hypothesis class. Let $\epsilon > 0$ and let $V_\epsilon = \mu(\{h : \varepsilon(h) \leq \epsilon\})$. Assume $\epsilon$ is large enough that $V_\epsilon > 0$. Let $\eta > 0$ and $\Delta \geq 0$ be such that $\Delta\eta \leq 1/2$. Then for any $\gamma \geq \ln(8/V_\epsilon)/\eta$,*

$$\Pr_{(x,y)\sim D}[y\ell(x) \leq 0] \leq 2(1 + (2/V_\epsilon)e^{-\eta\gamma})(\epsilon + \gamma),$$

*and*

$$\Pr_{(x,y)\sim D}[y\ell(x) \leq 2\Delta] \leq (1 + e^{2\Delta\eta})(1 + (2/V_\epsilon)e^{\eta(2\Delta - \gamma)})(\epsilon + \gamma)$$

$$\leq 4(1 + (2/V_\epsilon)e^{\eta(2\Delta - \gamma)})(\epsilon + \gamma).$$

The modification of Corollary 1 is immediate. In the discussion following Corollary 1, $\varepsilon(h^*)$ is replaced by $\epsilon$ as in Theorem 7.

Besides replacing $|\mathcal{H}|$ by $1/V_\epsilon$, the proof of Theorem 4 only needs to be modified by replacing all sums with integrals. Also, to upper bound $W_w$, we lower bound $Z$ by $V_\epsilon e^{-\eta\epsilon}$, a fact that follows immediately from the definition of $V_\epsilon$.

Generalizing the results of Section 6 to infinite class $\mathcal{H}$ seems harder and remains as an open problem for future research.

**8. Conclusions.** In this paper we present a new algorithm for prediction of binary functions using a weighted vote over all prediction rules within a class. We have shown when, and in what sense, this algorithm can perform better than the more common approach of choosing the prediction function which performs best on the training data.

While this algorithm is similar in spirit to a Bayesian prediction algorithm, there are at least two important differences.

The first difference is in the dependence of the posterior probability (before normalization) on the size of the training set $m$. In most Bayesian algorithms the expected value of the unnormalized posterior probability for any particular model $\Theta$ decreases at the rate $\exp(-c(\Theta)m)$, where $c(\Theta)$ is the expected value of the log probability of the data given the model. In our algorithm the rate of decrease is (approximately) $\exp(-c(\Theta)\sqrt{m})$. We choose this rate (Corollary 1) so that the variance of the empirical log-ratio is slowly decreasing, which results in an estimator whose stability improves as the size of the sample increases.

Second, the goal of our algorithm is to increase the stability of the prediction and not to optimize a Bayesian measure of risk. To that end, the only assumption regarding the data generation mechanism that we make in our analysis is that the data is generated in an IID fashion. To the best of our knowledge, all existing Bayesian analysis (other than on-line prediction methods) make the assumption that the data is generated by one of the models in the class over which the Bayesian averaging is performed. In this context it is worthwhile to mention recent work by Bousquet and Elisseeff [3] in which they show how improved generalization bounds can be proven for algorithms that are known to be stable. The main difference between that work and our work here is that we describe and analyze a specific averaging method that is guaranteed to be stable.

It was suggested that the main reason that our algorithm does not overfit has to do with the fact that we allow abstention, rather than with the averaging of many hypotheses. We believe that the most important property of our algorithm is the stability of the empirical log-ratio. Abstention is just one way of utilizing this stability. In other scenarios one may be better off using the log-ratio scores differently. For example, if the goal is to detect a rare type of instance within a large set, the correct method might be to sort all instances according to their log-ratio score and output the instances with the highest scores.

It is natural to think of the empirical log ratio as an estimate of the conditional probability of the label $y$ given the instance $x$. However, one should not take this intuition too far. The log ratio is a measure of the *model uncertainty* by which we mean the uncertainty in the identity of the best model which results from the finite size of the training set. It does *not* measure the uncertainty that is inherent in the true conditional distribution of $y$ given $x$. To realize this, consider a class with 100 rules in which one rule has a true error of 10%, while the true error of each of the other 99 rules is larger than 20%. Then with a training set with a few hundred examples the weight assigned to the best rule is likely to be larger than the total weight of all of the other rules. This in turn would imply that the log ratio would be very far from zero everywhere and our algorithm will always predict like the best rule and never abstain. Indeed, we can interpret the log-ratio values as an indication that we are certain which is the best rule in the class. This is quite independent of the fact that the best rule in the class has an error of 10%. To estimate this conditional probability we need to *calibrate* the predictions of our algorithm. One can devise various ways of performing this calibration. An interesting parameter-free calibration method has been recently suggested by Vovk [22].

Our work shares some ideas with the recent work by Shawe-Taylor and Williamson [20] and McAllester [16] on PAC–Bayesian analysis. The main common idea is that if many classification rules perform well, then their prediction can be trusted more than that of a single rule that is performing well. The main difference is that in our work we average over the predictions of the best rules and get a different prediction confidence for each test instance, while the PAC–Bayesian analysis uses the plurality of the good performers to improve the performance guarantees for a single classification rule that is chosen at random according to the posterior. Similar ideas were used in the analysis of large-margin classifiers.

Another connection worth mentioning here is to margin based classification methods such as SVMs [19, 21] and boosting [10, 18]. One intuition that explains why large margins are important regards the stability of the linear classifier. Large margins around the separating hyperplane imply that slight perturbations of the hyperplane will also classify the data correctly. In other words, it implies that a large set of similar linear classifiers have small training error. Suppose now that we used the averaging algorithm suggested in this paper where the set of classifiers that is used is the set of all linear classifiers. The fact that the set of close-to-optimal classifiers is large implies that the prediction where they all agree would be very confident. On the other hand, the region on which the algorithm will abstain is similar (but not identical) to the margin region. In other words, the behavior of our algorithm is, in fact, similar to that of large margin classifiers. However, there are two important differences. On the one hand, the averaging algorithm is much more general in that it can be applied to any set of classifiers, not just linear classifiers; neither does it depend on whether or not the data is separable, that is, perfectly classifiable by one of the rules in the class. On the other hand, our algorithm is extremely inefficient as compared to SVMs or AdaBoost as its application requires calculating the empirical error for each and every rule in the set.

**9. Future research.** We suggest two directions for future work, one regarding computational efficiency, the other regarding the choice of a prior distribution.

Consider first the computational issue. For most interesting hypothesis classes the task of finding the hypothesis that minimizes the training error is computationally intractable. Obviously, calculating the error of all of the hypotheses in the class is at least as hard as finding the best hypothesis and probably much harder. Does this mean that our algorithm cannot be used for practical learning problems? Not necessarily. Here are three approaches to solving the computational problem:

1. Sometimes the problem of learning a complex classification rule can be broken down into several problems of learning very simple rules. For example, Freund and Mason [9] show how to break down the problem of learning alternating decision trees (a class of rules which generalizes decision trees and boosted decision trees) into a sequence of simpler learning problems using boosting. Each of the simpler problems involves finding the best threshold rule in one dimension. These last problems are so simple the calculation can be done in time linear in the size of the training set. In this context our algorithm can be used directly and its use might significantly increase the robustness of the system as a whole.

2. In some cases a careful choice of the prior distribution over the hypotheses makes it possible to calculate the posterior average efficiently. For example, *conjugate priors* commonly used in Bayesian statistics are prior distribution which maintain their functional form as they are updated. A more interesting case which involves variable-length Markov models for sequences was studied by Willems, Shtarkov and Tjalkens [23] and extended by Helmbold and Schapire [12]. It might be possible to adapt these techniques to efficiently calculate the empirical log ratio for our algorithm.

3. In some cases the posterior distribution can be approximated by a single sharp peak around the best hypothesis. In such a case the empirical log ratio can be approximated using Laplace approximation method. This technique was used by Freund [8]. For an introduction to this type of approximation methods see the excellent book by de Bruijn [7].

4. Another approach to estimating the average vote over the empirically best hypothesis is to use random sampling. Suppose we are given a learning algorithm capable of finding a hypothesis with small training error. Our goal is to tweak the algorithm in a way that will randomly create a hypothesis whose performance is almost as good as the original untweaked hypothesis. Moreover, we want the distribution according to which the hypothesis is generated to be close to the distribution defined by our exponential weights.

There are several learning algorithms that sample hypotheses and average them. The best known of these so-called *ensemble* algorithms is Breiman's *bagging* algorithm [4, 5]. It might be that bagging is indeed an efficient randomized algorithm of the type suggested here. On the other hand, it might be possible to adapt the theory presented in this paper to give a rigorous analysis for the performance of bagging and other ensemble methods.

The second direction we suggest for future work is to consider the choice of the prior measure $\mu$ defined in Section 7. Clearly, the choice of measure has a large influence on the algorithm and on the upper bound given in Theorem 7.

Intuitively, we would like to maximize the probability measure of the set $V_\epsilon$. However, we need to define the measure $\mu$ *before* observing the training data, that is, before we know what $V_\epsilon$ is. One natural approach is to maximize the minimum over the measure of all possible sets $V_\epsilon$.

Consider first a case in which we have prior knowledge of the distribution over the instances, without the labels. In this case we can use the measure which places uniform weights over an $\epsilon$-net on the hypothesis class, as was suggested in Section 5. This will ensure that if the best hypothesis in the class has error $\varepsilon^*$, then the set $V_{\varepsilon^*+\epsilon}$ will have measure at least $1/N$ where $N$ is the size of the $\epsilon$-net. The disturbing thing about this choice for $\mu$ is that it depends on $\epsilon$. Possibly this disturbance can be cleared if one can use a limit distribution where $\epsilon \to 0$. Intuitively, such a limit measure will capture the detailed structure of the hypothesis space in a way similar to Jeffreys' prior in Bayesian analysis.

Assuming that this analysis can be carried through, one should return to the original problem in which the distribution over the instances is unknown. In this case we need to approximate the "ideal" algorithm by using the information about the instance distribution that we get from the training examples. Ultimately, we would like to find an averaging algorithm whose performance is close to the averaging algorithm that has this prior knowledge *and* that is efficiently computable.

## APPENDIX: PROOF OF LEMMA 4

First, let $\mathcal{K} = \{h_1, \ldots, h_N\}$, and let

$$F(\eta, \mathbf{x}) = \frac{1}{\eta} \ln \left( \sum_{i=1}^{N} e^{-\eta x_i} \right).$$

For any $\mathbf{x}$, by checking derivatives, it can be verified that the function $\eta \mapsto F(\eta, \mathbf{x})$ is nonincreasing, while the function $\eta \mapsto F(\eta, \mathbf{x}) - (\ln N)/\eta$ is nondecreasing. Therefore, if $0 < \eta_1 \leq \eta_2$, then for any $\mathbf{x} \in \mathbb{R}^N$,

$$(14) \qquad 0 \leq F(\eta_1, \mathbf{x}) - F(\eta_2, \mathbf{x}) \leq \left( \frac{1}{\eta_1} - \frac{1}{\eta_2} \right) \ln N.$$

Now let

$$\mathcal{E} = \left\{ \frac{4 \ln N}{i \lambda} : i = 1, \ldots, \left\lfloor \frac{4 \ln N}{\lambda} \right\rfloor \right\}.$$

We show next that for any $\eta \geq 1$, there exists $\hat{\eta} \in \mathcal{E}$ such that

$$\left| \frac{1}{\eta} - \frac{1}{\hat{\eta}} \right| \ln N \leq \frac{\lambda}{4}.$$

For if $\eta \geq 4(\ln N)/\lambda$, then let $\hat{\eta} = 4(\ln N)/\lambda$. Then

$$0 \leq \left( \frac{1}{\hat{\eta}} - \frac{1}{\eta} \right) \ln N \leq \frac{1}{\hat{\eta}} \ln N = \frac{\lambda}{4}.$$

Otherwise, if $1 \leq \eta \leq 4(\ln N)/\lambda$, then let $\hat{\eta} = 4(\ln N)/(i\lambda)$ be the smallest element of $\mathcal{E}$ that is no smaller than $\eta$. That is,

$$\frac{4 \ln N}{(i+1)\lambda} < \eta \leq \frac{4 \ln N}{i \lambda}.$$

Then

$$\begin{aligned}
0 \leq \left( \frac{1}{\eta} - \frac{1}{\hat{\eta}} \right) \ln N &= \left( \frac{1}{\eta} - \frac{i\lambda}{4 \ln N} \right) \ln N \\
&\leq \left( \frac{(i+1)\lambda}{4 \ln N} - \frac{i\lambda}{4 \ln N} \right) \ln N \\
&= \frac{\lambda}{4}.
\end{aligned}$$

Since $\hat{R}_\eta(\mathcal{K}) = F(\eta, \langle \hat{\varepsilon}(h_1), \ldots, \hat{\varepsilon}(h_N) \rangle)$, (14) and the argument above imply that for any $\eta \geq 1$, there exists $\hat{\eta} \in \mathcal{E}$ such that

$$|\hat{R}_\eta(\mathcal{K}) - \hat{R}_{\hat{\eta}}(\mathcal{K})| \leq \frac{\lambda}{4}$$

and so

$$|(\hat{R}_\eta(\mathcal{K}) - \mathrm{E}[\hat{R}_\eta(\mathcal{K})]) - (\hat{R}_{\hat{\eta}}(\mathcal{K}) - \mathrm{E}[\hat{R}_{\hat{\eta}}(\mathcal{K})])| \leq \frac{\lambda}{2}.$$

Thus,

$$\begin{aligned}
\mathrm{Pr}[\exists \eta \geq 1 &: |\hat{R}_\eta(\mathcal{K}) - \mathrm{E}[\hat{R}_\eta(\mathcal{K})] \mathcal{K}| \geq \lambda] \\
&\leq \mathrm{Pr}\left[ \exists \hat{\eta} \in \mathcal{E} : |\hat{R}_{\hat{\eta}}(\mathcal{K}) - \mathrm{E}[\hat{R}_{\hat{\eta}}(\mathcal{K})]| \geq \frac{\lambda}{2} \right] \\
&\leq 2|\mathcal{E}| e^{-\lambda^2 m/2},
\end{aligned}$$

where the second inequality uses the union bound combined with Lemma 1.

**Acknowledgments.** Most of the work done while the authors were working in AT&T Labs. Thanks to Michael Cameron-Jones and to the anonymous reviewers of this paper for helpful comments.

Y. Freund
Computer Science Department
Columbia University
New York, New York 10027
USA
e-mail: freund@cs.columbia.edu

Y. Mansour
Computer Science Department
Tel-Aviv University
Israel
e-mail: mansour@math.tau.ac.il

R. E. Schapire
Department of Computer Science
Princeton University
Princeton, New Jersey 08544
USA
e-mail: schapire@cs.princeton.edu